\documentclass{article}
\usepackage{amsmath,amssymb,amsthm}
\usepackage[pdftex]{graphicx}
\usepackage{fullpage}
\usepackage{cite}
\usepackage{palatino}
\usepackage[T1]{fontenc}
\usepackage{upquote}
\usepackage{color}

\usepackage{graphicx}
\usepackage[font=small,labelfont=bf]{caption}

\usepackage{subcaption}
\usepackage{amsmath,amsthm,amssymb}

\newcommand{\footremember}[2]{%
    \footnote{#2}
    \newcounter{#1}
    \setcounter{#1}{\value{footnote}}%
}

\newcommand{\R}{\mathbb{R}}

\newtheorem{thm}{Theorem}[section]
\newtheorem{lem}[thm]{Lemma}
\newtheorem{cor}[thm]{Corollary}
\newtheorem{prop}[thm]{Proposition}
\newtheorem{define}[thm]{Definition}
\newtheorem{rmk}[thm]{Remark}
\newtheorem{xmpl}[thm]{Example}

\renewcommand{\d}{\delta}

\newcommand{\lam}{\lambda}

\newcommand{\1}{\mathbf{1}}
\newcommand{\0}{\mathbf{0}}

\newcommand{\ip}[2]{\left\langle{{#1}},{{#2}}\right\rangle}

\renewcommand{\l}{\left}

\renewcommand{\r}{\right}
\newcommand{\be}{\begin{enumerate}}
\newcommand{\bi}{\begin{itemize}}
\newcommand{\ee}{\end{enumerate}}
\newcommand{\ei}{\end{itemize}}

\newcommand{\x}{{\mathbf {{x}}}}
\newcommand{\y}{{\mathbf {{y}}}}

\newcommand{\e}{{\mathbf {{e}}}}
\newcommand{\p}{{\mathbf {{p}}}}

\newcommand{\rv}{{\mathbf {{r}}}}

\newcommand{\diag}{\mathrm{diag}}

\newcommand{\In}{\mathrm{In}}

\title{A Locally Stable Equilibria Criterion \\
for the Generalized Lotka--Volterra 
}
\author{
Michael Richard Livesay \footremember{mrl}{Department of Mathematics, University of Illinois Urbana-Champaign, IL. }
}

\begin{document}

\maketitle

\begin{abstract}
The main result applies to non-degenerate cases of the generalized Lotka--Volterra model. 
A criterion is given that relates the stability of two fixed points with the associated Schur complement of there respective community matrices. 
\end{abstract}


%%%%%%%%%%%%%%%%%%%%%%%%%%%%%%%%%%%%%
%%%%%%%%%%%%%%%%%%%%%%%%%%%%%%%%%%%%%

\section{Introduction} \label{sec:model}

The Lotka--Volterra system is a biological model of quadratic equations which has been thoroughly studied; see \cite{hofbauer1998evolutionary}, \cite{hutson1983criterion}, \cite{hirsch1988systems}, \cite{zeeman1993hopf}, \cite{takeuchi1996global}, and \cite{zeeman1998three}. 
Given the number of species of the system the stable equivalence classes completely characterize the solutions of the system. 
For the case of $3$ species \cite{zeeman1993hopf} completely characterized the Lotka--Volterra into $33$ stable equivalence classes. 
A relationship is found between the Jacobian of fixed points, which show numerous restrictions in the stable equivalence classes in higher dimensions. 

Properties of the generalized Lotka--Volterra model are exploited yielding restrictions for stable equilibria. 
In Section \ref{sec:motiv} it is shown that the Jacobian at a fixed point on the boundary has two principal submatrices of interest; one correspond to the surviving species and the other to the extinct species. 
In Section \ref{sec:loc_stable} the relationships between those two principal submatrices are exploited to gain properties of the stability for the fixed point, see Theorem \ref{thm:gen_maxmin}. 
Theorem \ref{thm:plus_one} shows that if two equilibria are such that the sets corresponding to their surviving species are the same except that one has a single additional species more than the other, then at most only one of these equilibria is asymptotically stable.

The generalized Lotka--Volterra is defined for a given environment and set of $N$ interacting species. 
The environment applied to the species determines a real $N \times N$ {\em community matrix} $A = (A_{ij})$ and {\em population rate vector} $\rv \in \R^N$, where the population of species $i$ is given by $\x$. 
Then the generalized Lotka--Volterra is given by the dynamical system 
\begin{align}\label{eq:LV_ODE}
\frac{d}{dt}x_i = x_i \l( r_i - \sum_{j\in [N]} A_{ij} x_j \r), \forall i \in [N] \text{ and } x_i(0) > 0.
\end{align}
A fixed point equilibrium of (\ref{eq:LV_ODE}) must satisfy either $x_i = 0$ or $r_i = \sum_{j\in V} A_{ij} x_j,$ for each $i \in [N]$. 
If $A$ is nonsingular, then the solution to $\rv = A \x$ is unique.

The system of equations given by Eq \eqref{eq:LV_ODE} may also be written in vector form. 
The Schur product $\odot$ is the point-wise product 
$\odot\colon\R^N\times \R^N\to\R^N$ such that $(\x \odot \y)_i = x_iy_i$. 
Consider the following system 
\begin{align} \label{eq:A,r} 
\dot \x = \x \odot (\rv - A \x), \text{ with } \x(0) \in \R^N. 
\end{align}
The vector form of Eq \eqref{eq:LV_ODE} is Eq \eqref{eq:A,r} when $\x(0) \in \R_{>0}^N$. 
Equation Eq \eqref{eq:A,r} generalizes Eq \eqref{eq:LV_ODE} to allow some species to be extinct, or even have negative populations, which makes it possible to avoid claims with excessive special cases. 
Given $A, \rv$ the system given by Eq \eqref{eq:A,r} is often referred to as LV$(A, \rv)$ for convenience. 

If extinct, a species cannot reproduce, that is if $x_i(t) = 0$ for any $t >0$ then $x_i \equiv 0$. 
If not extinct, a species cannot die off in finite time, that is if $x_i(0)>0$ then $x_i(t)>0$ for all $t \in \R_{\ge 0}$. 
Nevertheless understanding Eq \eqref{eq:A,r} for the case when some species are extinct plays a role in the local behavior of Eq \eqref{eq:LV_ODE} where no species are extinct. 
This is clear by the form of the Jacobian of Eq \eqref{eq:A,r} at any fixed point. 
See Appendix \ref{apx:LV} for definitions of terms used.

%%%%%%%%%%%%%%%%%%%%%%%%%%%%%%%%%%%%%
%%%%%%%%%%%%%%%%%%%%%%%%%%%%%%%%%%%%%

\section{Projected Subsystems and Non-Degeneracy} \label{sec:notation}

The notation is somewhat terse and is motivated by the frequent use of principal submatrices and Schur complements. 
Assume $N \times N$ matrix $A$ and $\p \in \R^N$ is given. 
Removing a given set of columns and their associated rows from the community matrix, as well as the same rows of the growth vector, gives a subsystem of LV$(A, \rv)$. 
Consider the non-empty ordered set $S \subset [N]$, with order $\le$, and let $S_i$ denote the $i^{th}$ largest element of $S$. 
Denote $P^S$ as the $N \times |S|$ matrix defined so that 
\begin{equation}
P^S_{i,j} = \d_{i,S_j}. 
\end{equation}
Denote the associated projection 
\begin{equation}
P_S := P^S \left(P^S \right)^*.
\end{equation}
For example, if $S = \{1,3,4\}$ with $N=4$, then
\begin{align*} 
P^{\{1,3,4\}} = 
\begin{pmatrix} 
1 & 0 & 0 \\
0 & 0 & 0 \\
0 & 1 & 0 \\
0 & 0 & 1 \\
\end{pmatrix}, 
P_{\{1,3,4\}} = 
\begin{pmatrix} 
1 & 0 & 0 & 0 \\
0 & 0 & 0 & 0 \\
0 & 0 & 1 & 0 \\
0 & 0 & 0 & 1 \\
\end{pmatrix}.
\end{align*}
In general $P_S = \diag( \sum_{i \in S} \e_i)\le I$. 
Rearranging the columns and rows of these matrices allows tools such as determinant products and Schur complements to be applied in more generic ways. 
Consider a second non-empty ordered set $S,T \subseteq [N]$. 
Define 
\begin{equation}
\begin{aligned} 
P^{S,T} := 
\begin{pmatrix} 
|	 	&  	& |	 		& |	 		&  	& | \\
\e_{S_1} 	& ... 	& \e_{S_{|S|}} 	& \e_{T_{1}} 	& ... 	& \e_{T_{|T|}} \\
|	 	&  	& |	 		& |	 		&  	& | \\
\end{pmatrix}, 
\text{ for disjoint $S$ and $T$.}
\end{aligned}
\end{equation}
Notice that $P^{S,T} \left(P^{S,T}\right)^* = P_{S \cup T} = P^{S \cup T} \left(P^{S \cup T}\right)^* $. 

Denote the pseudo inverse of a matrix, $M$, by $M^\dagger$. 
Suppose that $S,T \subseteq [N]$ are non-empty, $i \in {(S \cup T)}^c$ and that $\left(P^S\right)^* A P^S, \left(P^T\right)^* A P^T$ are invertible principal submatrix of $A$, then denote the following well-defined terms: 
\begin{equation} 
\begin{aligned} %\label{proj_vars}
A^S &:=(P^S)^* A P^S  			&A_S &:= P_S A P_S  \\
\rv^S &:=(P^S)^* \rv  				&\rv_S &:= P_S \rv  \\
\p^S &:= (A^S)^{-1} \rv^S 			& \p_S &:= A_S^{\dagger} \rv_S \\
%D^S(\x) &:= \diag \l(\l(P^S\r)^* \x\r) 
\end{aligned}
\end{equation}
\begin{equation}
\begin{aligned}
A^{S,T} &:= \left(P^{S,T}\right)^* A P^{S,T} = 
\begin{pmatrix}
(P^S)^* A P^S 	&	(P^S)^* A P^T \\
(P^T)^* A P^S	& 	(P^T)^* A P^T \\
\end{pmatrix},
&\text{ for disjoint $S$ and $T$.} \\
A^{S,i} &:= A^{S,\{i \}} =
\begin{pmatrix}
(P^S)^* A P^S 	&	(P^S)^* A \e_i \\
\e_i^* A P^S	& 	\e_i^*A \e_i \\
\end{pmatrix},
&\text{ for $i \notin S$.} \\
\end{aligned}
\end{equation}
%The subscript is always dropped in the case where $S = [N]$. ????
%, for example $D:= \diag(\p)$. 
%If $i > j$ for each $j \in S$ then the Schur complement 
%$A^{S \cup \{i\}}/ A^S = A_{ii} - \e_i^* A A_S^{\dagger}  A \e_i$. 
Notice that 
\begin{align*}
\p_S &= P^S \p^S, \\
A_S^{\dagger} &= P^S((P^S)^* A P^S)^{-1} (P^S)^*, \\
A^{S,i} / A^S &= A_{ii} - \e_i^* A A_S^{\dagger}  A \e_i, \text{ for $i \notin S$.}  
\end{align*}
If ever $A^S$ is not invertible, one could define $\p^S$ as $(A^S)^{\dagger} \rv^S$ instead. 
%Hence given $A$ and either $\rv$ or $\p$, the other is inferred as well as $\rv_S$ and $\p_S$ for any non-empty $S \subseteq [N]$. 
Notice that $A^{S,i}$ is simply one row and one column operation applied to $A^{S \cup \{ i \}}$. 
The vector $\p_S$ may not always be in the domain of Eq \eqref{eq:LV_ODE}, but it is always in the domain of Eq \eqref{eq:A,r}. 
%Even if $\p_S \ge \0$ it may not be stable, for example if $\diag(\p_S) A_S$ has a negative eigenvalue. 

\begin{rmk} \label{rm:sub_dyn}
Suppose that $S$ is given.
If $\x = P_S \x$ then 
\begin{align} \label{eq:sub_dyn}
\dot{\x} = \x \odot (\rv - A \x) = P_S \x \odot (\rv - A P_S \x) = P_S \x \odot (P_S \rv - P_S A P_S \x). 
\end{align}
Thus $\p_S = (P_S A P_S)^{\dagger}\rv$ is a fixed point for LV$(A, \rv)$. 
\end{rmk}

%%%%%%%%%%%%%%%%%%%%%%%%%%%%%%%%%%%%%
%%%%%%%%%%%%%%%%%%%%%%%%%%%%%%%%%%%%%

Remark \ref{rm:sub_dyn} shows how to find all of the fixed points for LV$(A, \rv)$. 
To illustrate Remark \ref{rm:sub_dyn} consider the following example. 

\begin{xmpl} \label{xpl:nonsingular_reduced}
Let 
\begin{align*}
A = 
\begin{pmatrix}
2 & -1 & 2 \\
1 & 2 & 1 \\
-3 & 1 & 1 \\
\end{pmatrix}
\text{ and }
\rv = 
\begin{pmatrix}
2 \\
4 \\
1 \\
\end{pmatrix}. 
\end{align*}
The fixed point of Eq \eqref{eq:A,r} associated with the surviving species given by $S=\{1,2\}$. 
\begin{align*}
P_S A P_S = 
\begin{pmatrix}
2 & -1 & 0 \\
1 & 2 & 0 \\
0 & 0 & 0 \\
\end{pmatrix},
P_S \rv = 
\begin{pmatrix}
2 \\
4 \\
0 \\
\end{pmatrix},
\text{ and } \y = (P_S A P_S)^{\dagger} \rv = 
\begin{pmatrix}
8/5 \\
6/5 \\
0 \\
\end{pmatrix}. 
\end{align*}
Then $\y$ is a fixed point of Eq \eqref{eq:A,r} because 
\begin{align*}
\begin{pmatrix}
8/5 \\
6/5 \\
0 \\
\end{pmatrix}
\odot \l(\rv - A 
\begin{pmatrix}
8/5 \\
6/5 \\
0 \\
\end{pmatrix}\r)=
\begin{pmatrix}
8/5 \\
6/5 \\
0 \\
\end{pmatrix}
\odot 
\begin{pmatrix}
0 \\
0 \\
23/5 \\
\end{pmatrix}=\0. 
\end{align*}
\end{xmpl} 

Consider Eq \eqref{eq:A,r} and some set $S \subseteq [N]$. 
The principal submatrix associated with $S$ is invertible if and only if the solution $\x$ to $A_S \x = P_S \rv$ subject to $\x = P_S \x$ is unique. 
Thus the number of fixed points to Eq \eqref{eq:A,r} is at most $2^N$ if every principal submatrix of $A$ is invertible, and otherwise the number of fixed points is uncountable. 
Due to such differences it is important to know whether or not all the principal submatrices are nonsingular; both cases are interesting and should be studied separately. 
The case where all principal submatrices are invertible is the more generic case, and the other case is not covered here. 
Even a single singular principal submatrix is not natural without assuming additional structure for the system. 
Indeed, if a small amount of noise is added along the diagonal of $A$ then any principal submatrix of $A$ would be nonsingular with probability $1$. 
To that end, consider the following definitions, see \cite{murty1972number, murty1988linear}. 

%%%%%%%%%%

\begin{define} \label{def:nondeg_A}
A matrix is {\em non-degenerate} if all its principal submatrices are nonsingular. 
\end{define}

\begin{define} \label{def:nondeg_r}
Given the matrix $A$, the vector $\rv$ is non-degenerate with respect to $A$ 
if $\0 = \x \odot \l( \rv - A \x \r)$ implies that for each $i \in [N]$, exactly one of $x_i = 0$ or $r_i - \sum_{j\in [N]} A_{ij} x_j = 0$. 
\end{define}

Together these non-degeneracies are very powerful. 
For simplicity consider the following definition. 

\begin{define} \label{def:nondeg_pair} 
Given the matrix $A$ and vector $\rv$, the pair $(A,\rv)$ is called {\em non-degenerate} if $A$ is non-degenerate and if $\rv$ is non-degenerate with respect to $A$. 
\end{define}

This non-degenerate case implies that Eq \eqref{eq:A,r} has exactly $2^N$ distinct fixed points, as demonstrated in the works of \cite{murty1972number, murty1988linear}, as laid out in Theorem \ref{thm:nondeg_pair}. 

%%%%%%%%%%

\begin{thm} \label{thm:nondeg_pair}
Suppose $A$ is an $N \times N$ matrix and $\rv \in \R^N$. 
Then $(A, \rv)$ is non-degenerate if and only if there are exactly $2^N$ distinct solutions to 
$\x \odot (\rv - A \x) = \0$. 
\end{thm}

%%%%%%%%%%%%%%%%%%%%%%%%%%%%%%%%%%%%%%%%%%%%%

\iffalse
\begin{proof}

Fix $S \subseteq [N]$. 
Assume that $(A, \rv)$ is non-degenerate. 
Since $A$ is non-degenerate, for any $S \in [N]$ there is only one solution to $P_S \rv - A_S \x = \0$ for $\x = P_S \x$. 
Also $\x \odot (\rv - A \x) = P_S \x \odot (\rv - A P_S\x) = \x \odot (P_S \rv - A_S \x) = \0$. 
This implies that there are at most $2^N$ solutions to $\x \odot (\rv - A \x) = \0$. 
Consider $S \neq T \subset [N]$. 
If $\x$ solves $P_S (\rv - A \x) = \0$ with $\x = P_S \x$ then it cannot solve $P_T (\rv - A \x) = \0$ with $\x = P_T \x$, because $(A, \rv)$ is non-degenerate. 
Thus there are $2^N$ solutions to $\x \odot (\rv - A \x) = \0$.

%%%%%%%%%%%%%%%%%%%%%%%%%%%%%%%%%%%%%%%%%%%%%

For the other direction first suppose $A$ is degenerate. 
There exists an $S \in [N]$ so that $A^S$ is singular and so there are infinitely many solutions to $P^S \rv - A^S \x = \0$ for $\x \in \R^{|S|}$ or no solutions at all. 
Thus there are less than $2^N$ solutions to $\x \odot (\rv - A \x) = \0$ or infinitely many. 

Now suppose that $A$ is non-degenerate and $(A, \rv)$ is degenerate. 
Recall that this implies that there are at most $2^N$ solutions. 
Since $(A, \rv)$ is degenerate there exist $S \subset T \in [N]$ and a solution $\x$ to $P_T \rv - A_T \x = \0$ and $P_S \rv - A_S \x = \0$ and $x_i =0$ for $i \in T \backslash S$. 
Thus there are at most $2^N -1$ solutions. 
The conclusion follows. 
\end{proof}
\fi
%%%%%%%%%%%%%%%%%%%%%%%%%%%%%%%%%%%%%%%%%%%%%

Another interesting implication is that if $(A,\rv)$ is non-degenerate then the Jacobian at each equilibria is nonsingular, see Theorem \ref{thm:FT_stable_crit}. 
Working in the paradigm of non-degenerate pairs $(A, \rv)$ simplifies the analysis. 
See Appendix \ref{apx:props} for basic properties of the terms denoted in this section.

\section{Jacobian Criterion for Locally Stable Equilibria} \label{sec:motiv}

%%%

%Global stability results for Lotka--Volterra dynamics have been thoroughly studied. 
%Certain properties of global stability of the reduced system can be extended to give restrictions of locally stable equilibria, due to properties of the form of the Jacobian matrix. 

%%%

For the general Lotka--Volterra system, the environment determines the community matrix and growth rate vector for any set. 
Suppose that for this set, the Lotka--Volterra model has a globally stable equilibrium in the positive orthant. 
This globally stable equilibria often tells little about its local stability in the system after introducing additional species, as seen in Examples \ref{xpl:not_stable_subs} and \ref{xpl:bdry_fp}. 
Further analysis of locally stable equilibria is explored in this section. 
In this section a criterion for fixed points to be stable for the Lotka--Volterra system is discussed. 
Assume that there are $N$ species throughout.

%%%%%%%%%%%%%%%%%%%%%%%%%%%%%%%%%%%%%%%%%%%%%

\begin{xmpl} \label{xpl:not_stable_subs}
Consider an environment with $3$ species and community matrix and growth rate vector 
\begin{align*}
A =
\begin{pmatrix}
5	& 	-2 	& 	-2  \\
-2	& 	5 	& 	-2  \\
-2	& 	-2 	& 	5  \\
\end{pmatrix}
\text{ and }
\rv = 
\begin{pmatrix}
1 \\
1 \\ 
1 \\
\end{pmatrix}.
\end{align*}
Since $A$ is symmetric $V(\x)=\x^* A \x$ is a Lyapunov function for this system if and only if $A$ is positive definite. 
Thus, it is clear that any reduced system with less than $3$ species has a global equilibrium, since $A_S > 0$ when $|S|<3$. 
Yet only the equilibrium $\1$ is stable in the full system as it is globally stable by Lyapunov's theorem, hence $\p_S$ is not stable.
Indeed, $A(\p - \p_{S}) \ge 0$ for $|S|<3$, thus $\p_{S}$ is not even a {\em saturated} fixed point, see Fig. \ref{fig:first_xpl}.

\begin{figure}\centering
\begin{subfigure}{.42\textwidth}\centering
\includegraphics[width=1\linewidth]{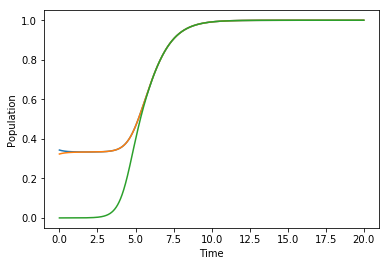}
\end{subfigure}%
\begin{subfigure}{.42\textwidth}\centering
\includegraphics[width=1\linewidth]{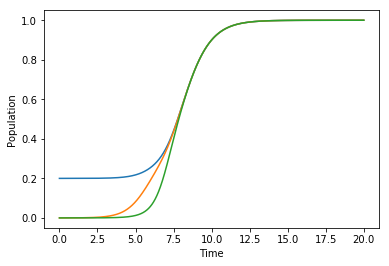}
\end{subfigure}%
\caption{
Plotted are two instances with initial conditions near two different fixed points, showing each species' population converges to $1$. 
On the left $|S|=2$ and on the right $|S|=1$, in Ex. \ref{xpl:not_stable_subs}. 
} 
\label{fig:first_xpl}
\end{figure}

\end{xmpl} 

In Ex \ref{xpl:bdry_fp} it is shown that even if $ \p^S$ is stable in the reduced system LV$(A^S, \rv^S)$ and $\p_T$ is not stable for every $T$ where $S \subset T$, it is not necessarily true that $\p_S$ is stable in LV$(A, \rv)$. 
The converse is also false, that is there exists examples so that $S \subset T$ and both $\p_S$ and $\p_T$ are stable, see Example \ref{xpl:counter}. 

\begin{xmpl} \label{xpl:bdry_fp}
Consider an environment and $4$ species. 
Assume species $4$ is extinct. 
Assume that the community matrix and growth rate vector of the remaining three species are given as

\begin{align*}
A_{\{1,2,3\}} =
\begin{pmatrix}
1	& 	.5 	& 	0 	& 	0  \\
.5	& 	1 	& 	0 	& 	0  \\ 
0	& 	0 	& 	1 	& 	0  \\
0	& 	0 	& 	0 	& 	0  \\
\end{pmatrix}
\text{ and }
\rv_{\{1,2,3\}} = 
\begin{pmatrix}
1.3 \\
1.3 \\ 
0.8 \\
0 \\
\end{pmatrix}.
\end{align*}
The system, LV$\l(A^{\{1,2,3\}}, \rv^{\{1,2,3\}}\r)$ has a globally stable equilibrium at 
$\p^{\{1,2,3\}} = 
\begin{pmatrix}
0.8666... \\
0.8666... \\ 
0.8 \\
\end{pmatrix}$. 

\begin{figure}\centering
\begin{subfigure}{.42\textwidth}\centering
\includegraphics[width=1\linewidth]{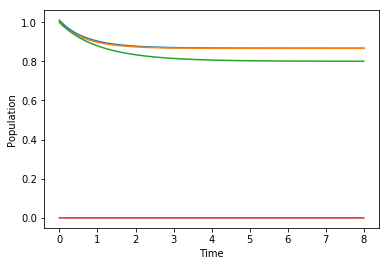}
\end{subfigure}%
\caption{
The fourth species is left out of the environment, and the other $3$ species prosper, in Ex. \ref{xpl:bdry_fp}. 
} 
\label{fig:fourth_extinct}
\end{figure}

Suppose the fourth species is reintroduce back into the environment. 
Then depending on how it interacts with the other species, completely different dynamics can be seen. 
It should be noted that the new interactions are completely independent of the currently observed parameters. 
Suppose that the full community matrix and growth rate vector are 
\begin{align*}
A =
\begin{pmatrix}
1	& 	.5 	& 	0 	& 	2 \\
.5	& 	1 	& 	0 	& 	2 \\ 
0	& 	0 	& 	1 	& 	2 \\
2	& 	2 	& 	2 	& 	1 \\
\end{pmatrix}
\text{ and }
\rv = 
\begin{pmatrix}
1.3 \\
1.3 \\ 
0.8 \\
5.9 \\
\end{pmatrix}, 
\text{ hence }
\p = 
\begin{pmatrix}
1 \\
1 \\ 
1 \\
-0.1 \\
\end{pmatrix}, 
\end{align*}
The Jacobian at $\p$, given by Eq \eqref{eq:Jacobian}, has a positive eigenvalue and thus $\p$ is not stable. 
The points $
\begin{pmatrix}
0.8666... \\
0.8666... \\ 
0.8 \\
\end{pmatrix}$ and $(5.9)$ are asymptotically stable in the subsystems with existing species $\{1,2,3\}$ and $\{4\}$ respectively. 
Only the fixed point $\p_{\{4\}}$ is saturated in LV$(A,\rv)$. 
In fact $\p_{\{1,2,3\}}$ is not even a saturated equilibrium since $\e_4^* A(\p - \p_{\{1,2,3\}}) > 0$, see Fig. \ref{fig:fourth}. 
Having even a small initial value for $x_4$ makes it impossible for the point $\p_{\{1,2,3\}}$ to be attracting.

\begin{figure}\centering
\begin{subfigure}{.42\textwidth}\centering
\includegraphics[width=1\linewidth]{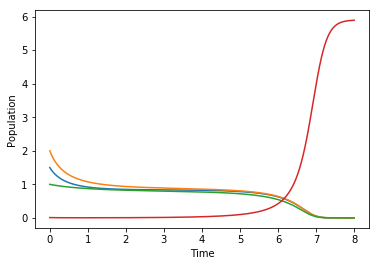}
\end{subfigure}%
\begin{subfigure}{.42\textwidth}\centering
\includegraphics[width=1\linewidth]{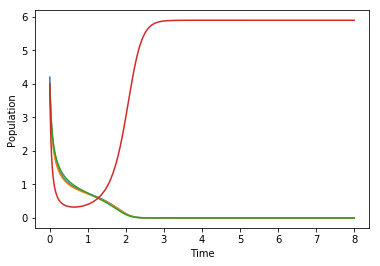}
\end{subfigure}%
\caption{
To the left, the dynamics appears to converge to $\p_{\{1,2,3\}}$ at first, but over time the population of $x_4$ eventually prevails. 
To the right, the initial populations of each species is larger. 
The value of $x_4$ starts off diminishing even at a greater rate than the other species, but unavoidably, the population $x_4$ persists in defiance of the other three, in Ex. \ref{xpl:bdry_fp}. 
} 
\label{fig:fourth}
\end{figure}

\end{xmpl}

%%%%%%%%%%%%%%%%%%%%%%%%%%%%%%%%%%%%%%%%%%%%%

Example \ref{xpl:bdry_fp} considers stability in a reduced system and shows that even if there is no strictly larger system with a stable point in its interior that the reduced system's stable point may not be stable in the full system.  
However, the fixed point in the full system was not a saturated fixed point. 

%%%%%%%%%%%%%%%%%%%%%%%%%%%%%%%%%%%%%%%%%%%%%

The fixed point $\p_S$ is stable if and only if $\p_S$ is a saturated fixed point in LV$(A, \rv)$ and $\p^S$ is stable in LV$(A^S, \rv^S)$ corresponding to $S$, see Theorem \ref{thm:FT_stable_crit}. 
To that end let $g_i(\x) = x_i (r_i - \sum_{j=1}^N a_{ij} x_j)$, then 
\begin{equation}
\begin{aligned}
\frac{\partial}{\partial x_j} g_i(\x) &= 
\begin{cases}
-a_{ij} x_i,	&	\mbox{for } i \neq j \\
r_i - a_{ii} x_i - \sum_{k = 1}^N a_{ik} x_k,	&	\mbox{for } i = j  \\
\end{cases}\\
\frac{\partial}{\partial x_j} g_i(\x) &= (A(\p - \x))_i \delta_{ij} - a_{ij} x_i 
\end{aligned}
\end{equation}
In general the Jacobian, denoted $J$, at $\x$ is 
\begin{align} \label{eq:Jacobian}
J(\x) = \diag(A(\p - \x)) - \diag(\x) A. 
\end{align}

By the form of the Jacobian, stability requires that $A(\p - \x) \le \0$. 
It is mentioned in \cite{takeuchi1996global} that the saturated fixed points of the Lotka--Volterra are the solutions to the linear complementarity problem; see Proposition \ref{prop:LV_LCP}. 
A matrix is called stable if all of its eigenvalues are negative.

\begin{thm} \label{thm:FT_stable_crit}
Suppose $(A,\rv)$ is non-degenerate, given an $N \times N$ matrix, $A$, and $\rv \in \R^N$. 
For any non-empty $S \subseteq [N]$ the Jacobian $J(\p_S)$ is nonsingular. 
Moreover, the point $\p_S \in \R^N$ is asymptotically stable for the system LV$(A,\rv)$
iff $A(\p - \p_S) \le \0$ and the matrix $-\diag(\p^S) A^S$ is stable. 
\end{thm}

\begin{proof}
Without loss of generality suppose that there is a $k < N$ so that $S = [k]$, that is $A = A^{S, S^c}$. 
Notice that $\diag(A(\p - \p_S)) = P_{S^c} \diag(A (\p - \p_S))$ and that $\diag(\p_S) A = P_S \diag(\p_S) A$. 
Thus 
\begin{equation} \label{eq:Jacob_p}
\begin{aligned}
J(\p_S) 
&= \diag(A(\p - \p_S)) - \diag(\p_S) A \\ 
&=
\begin{pmatrix}
0_{|S| \times |S|}			&	0_{|S| \times |S^c|} \\
0_{|S^c| \times |S|}				&	\diag( A(\p - \p_S))^{S^c}  \\
\end{pmatrix} -
\begin{pmatrix}
\diag(\p^S) A^S			&	B \\
0_{|S^c| \times |S|}				&	0_{|S^c| \times |S^c|}  \\
\end{pmatrix}, \\
&=
\begin{pmatrix}
-\diag(\p^S) A^S			&	-B \\
0_{|S^c| \times |S|}				&	\diag( A(\p - \p_S))^{S^c}  \\
\end{pmatrix}, \text{ for some matrix } B. 
\end{aligned}
\end{equation}
Note the entries of $(P^{S^c})^* A(\p - \p_S)$ are called {\em transversal eigenvalues}. 
By Proposition \ref{prop:Schur_det_id} 
\begin{align*}
\det\left(J(\p_S) - \lam I\right) &= 
\det \left(-\diag(\p^S) A^S - \lam I_{|S| \times |S|} \right) \cdot
\det \left( \diag( A(\p - \p_S))^{S^c} - \lam I_{|S^c| \times |S^c|} \right). 
\end{align*}
The eigenvalues of $J(\p_S)$ are exactly the accumulation of eigenvalues of $-\diag(\p^S) A^S$ and the entries of $(P^{S^c})^* A(\p - \p_S)$. 
Since $\det \l(-\diag(\p^S) A^S \r) = \det \l(-\diag(\p^S) \r) \det\l(A^S \r)$, the eigenvalues of $-\diag(\p^S) A^S$ are non-zero. 
Since each input of $(P^{S^c})^* A(\p - \p_S)$ is nonzero, $\det \left( \diag( A(\p - \p_S))^{S^c} \right) \neq 0$, and so $J(\p_S)$ is nonsingular. 

Finally since $A(\p - \p_S) \le \0$ if and only if  $\diag( A(\p - \p_S))^{S^c} < 0$. 
We have that $J(\p_S)$ is stable if and only if $A(\p - \p_S) \le \0$ and the matrix $-\diag(\p^S) A^S$ is stable.

\iffalse
Suppose $S = [N-1]$. 
Then 
\begin{align*}
\e_N^* A(\p - \p_S) &= \e_N^* A (\p - (P_S A P_S)^{\dagger} \rv ) \\
&= \e_N^* (A - A (P_S A P_S)^{\dagger} A) \p \\
&= \e_N^* \l( 
\begin{pmatrix}
A^S 	& A \e_N \\
\e_N^*A 		& A_{NN} \\
\end{pmatrix}
 - 
\begin{pmatrix}
A^S 	& A \e_N \\
\e_N^*A 		& A_{NN} \\
\end{pmatrix} 
\begin{pmatrix}
A^S 	& \0 \\
\0^T 		& 0 \\
\end{pmatrix} 
\begin{pmatrix}
A^S 	& A \e_N \\
\e_N^*A 		& A_{NN} \\
\end{pmatrix} 
 \r) \p \\
&= \e_N^*
\begin{pmatrix}
0		 	& \0 \\
\0^T 			& A_{NN} - \e_N^* A A_S A \e_N \\
\end{pmatrix}
\p \\
&= (A_{ii} - \e_N^* A A_S^{\dagger}  A \e_N) (\e_N^* \p) \\ 
&\neq 0, \mbox{ since $(A,\rv)$ is non-degenerate $\e_N^* \p \neq0$ and so Proposition \ref{prop:comp_nondeg} applies.}
\end{align*}
Thus for any non-empty $S \subset [N]$ and $i \notin S$, $\e_i^* A(\p - \p_S) \neq 0$. 
\fi

\end{proof}

\begin{rmk} \label{rmk:crit}
There is another way to word Theorem \ref{thm:FT_stable_crit} for those familiar with LCPs. 
Suppose the pair $(A,\rv)$ is non-degenerate. 
Fix a non-empty $S \subset [N]$, then the Jacobian $J(\p_S)$ is nonsingular. 
Furthermore, $\p_S$ is a stable equilibrium of LV$(A,\rv)$ iff $\p_S$ solves LCP$(A,\rv)$, see definition \ref{def:LCP}, and $\p^S$ is stable in LV$(A^S,\rv^S)$. 
\end{rmk}

The assumption that $S$ was non-empty in Theorem \ref{thm:FT_stable_crit} was due to there being no meaning assigned to $A^{\{\}}$. 
The special case where $S$ is empty is simple to handle. 
Suppose $(A,\rv)$ is a non-degenerate. 
The Jacobian at $\0$ is $J(\0) = \diag \left( A \p \right)$, so all species go extinct, near $\0$, if and only if $\rv = A \p < \0$. 

If one was to attempt to generalize Theorem \ref{thm:FT_stable_crit} to the case $(A,\rv)$ is a not non-degenerate, they would have to additionally consider the nature of the uncountable fixed points. 
This extension is not covered here.

The next example shows that if $(A, \rv)$ is not non-degenerate, as in Theorem \ref{thm:FT_stable_crit}, then the Jacobian at $\p_S$ could be singular. 

\begin{xmpl} \label{xpl:singular_jac}
Consider a set $S$ such that  $\p_S$ is a saturated fixed point and $\p^S$ is stable in its corresponding reduced system. 
It is possible that the Jacobian at $\p_S$ is singular, giving cause to assume $(A, \rv)$ is non-degenerate. 

Consider an environment with $3$ species and with the non-degenerate community matrix and growth rate vector given as 
\begin{align*}
A =
\begin{pmatrix}
1	& 	0 	& 	3  \\
-2	& 	1 	& 	0  \\
2	& 	0 	& 	1  \\
\end{pmatrix}
\text{ and }
\rv = 
\begin{pmatrix}
1 \\
-2 \\ 
1 \\
\end{pmatrix}.
\end{align*}
Then $\p_{\{1\}} = 
\begin{pmatrix}
1 \\ 
0 \\ 
0 \\ 
\end{pmatrix}$ and
$\p = 
\begin{pmatrix}
0.4 \\ 
-1.2 \\ 
0.2 \\ 
\end{pmatrix}$.
So $A(\p - \p_{\{1\}}) = 
\begin{pmatrix}
0 \\ 
0 \\ 
-1 \\ 
\end{pmatrix} \le \0$, thus $\p_{\{1\}}$ is a saturated equilibrium. 
Clearly $\p^{\{1\}}=1$ is stable in the reduced system because $1(r_1-A_{\{1\}}1)=0$, and $\p_{\{1\}}= \p_{\{1,2\}}= (1,0,0)^T$. 
Yet, the Jacobian at $\p_{\{1\}}$ is singular since 
$\diag(A(\p - \p_{\{1\}})) - \diag(\p_{\{1\}}) A= 
\begin{pmatrix}
0	& 	0 	& 	0  \\
0	& 	0 	& 	0  \\
0	& 	0 	& 	-1  \\
\end{pmatrix}
-
\begin{pmatrix}
1	& 	0 	& 	0  \\
0	& 	0 	& 	0  \\
0	& 	0 	& 	0  \\
\end{pmatrix} A
=
\begin{pmatrix}
-1	& 	0 	& 	-3  \\
0	& 	0 	& 	0  \\
0	& 	0 	& 	-1  \\
\end{pmatrix}
$; see Eq \eqref{eq:Jacobian}. 
\end{xmpl} 

As shown in Theorem \ref{thm:FT_stable_crit}, assuming that $(A, \rv)$ is non-degenerate rectifies the issue of singular Jacobians at equilibria. 
Section \ref{sec:loc_stable} connects the properties of the populations of the living with the dead in the Jacobian at an equilibrium.

%\pagebreak

%%%%%%%%%%%%%%%%%%%%%%%%%%%%%%%%%%%%%
%%%%%%%%%%%%%%%%%%%%%%%%%%%%%%%%%%%%%

\section{Schur Complement and Local Stability} \label{sec:loc_stable}

The main results of this section are given depends on Proposition's \ref{prop:comp_form} and \ref{prop:stable_form}, which are proven in Appendix \ref{apx:props}. 
Suppose that the pair $(A, \rv)$ is non-degenerate with $N$ species and $S = [k] \subset [N]$. 
Then Proposition \ref{prop:comp_form} gives that the transversal eigenvalues with respect to $S$ are 
\begin{align} \label{eq:transversal_eval}
\begin{aligned}
\l(P^{S^c}\r)^* A(\p - \p_S) 
= \l( A / A^S \r) \l(P^{S^c}\r)^* \p. 
%= \l( A / A^S \r) \diag(Q_{S^c}^* \p) \1. 
\end{aligned}
\end{align}
While Proposition \ref{prop:stable_form} gives that 
\begin{align} \label{eq:main_eval}
\begin{aligned}
\l( \diag(\p) A \r) / \l( \diag(\p) A \r)^S = \diag\l(\l(P^{S^c}\r)^* \p \r) \l( A / A^S \r). 
\end{aligned}
\end{align}
So by Eq \eqref{eq:Jacob_p}, \eqref{eq:transversal_eval}, and \eqref{eq:main_eval} we have that 
\begin{align} \label{eq:eloq}
\begin{aligned}
%J(\p_S)^{S} &= -\diag(\p^{S}) A^S, \text{ and } \\
%J(\p_S)^{S^c} &= \diag \l(\l( A / A^S \r) \diag(\p)^{S^c} \1 \r), \text{ and } \\
J(\p_S)^{S^c} &= \diag \l(\l( A / A^S \r) \l(P^{S^c}\r)^* \p  \r) \\
J(\p)/ J(\p)^S &= \diag\l(\l(P^{S^c}\r)^* \p \r) \l( A / A^S \r). 
\end{aligned}
\end{align}
This means that the Jacobian at $\p_S$ relates to the Jacobian at $\p$. 
Specifically, this shows that there is a strong relationship between the transversal eigenvalues (i.e. the diagonal elements of $J(\p_S)^{S^c}$) corresponding to $S$ and the Schur complement (i.e. $J(\p)/ J(\p)^S$) of the Jacobian at $\p$ over its principal submatrix corresponding to $S$. 
This is not obvious and can be useful as seen in Example \ref{ex:last_ex}. % make example SEE EXAMPLES
In Theorem \ref{thm:gen_maxmin} this relationship is explored further.

%%%%%%%%%%%%%%%%%%%%%%%%%%%%%%%%%%%%%

\begin{xmpl} \label{ex:last_ex}
Consider the non-degenerate $(A, \rv),$ \\
$A =
\begin{pmatrix}
2	& 	-2 	& 	4 \\
0	& 	2 	& 	-1  \\
0	& 	-1 	& 	2  \\
\end{pmatrix},$
$\rv = 
\begin{pmatrix}
2 \\
2 \\
2 \\
\end{pmatrix},$ so $\p = (-1,2,2)^T$. 
Let us analyze if either $\p_{\{1\}}$ or $\p_{\{2,3\}}$ are stable. 

Let $S = \{1\}$, 
then 
$A/A^S = 
\begin{pmatrix}
2 	& 	-1  \\
-1 	& 	2  \\
\end{pmatrix}$
and so 
$A/A^S P^{S^c} \p = 
\begin{pmatrix}
2  \\
2  \\
\end{pmatrix}$, thus $\p_{\{1\}} = (1,0,0)^T$ is unstable by Proposition \ref{prop:comp_form}. 

Let $T = \{2,3\}$, then $A/A^T = 2$ and so $A/A^T P^{T^c} \p = -2$. 
Also, $\p_T = (0,2,2)^T$ and so $-\diag(\p^T) A^T$ is stable. 
\end{xmpl}

In Example \ref{ex:last_ex}, $\p_{\{2,3\}}$ is stable, and thus $\p_{\{2\}}, \p_{\{3\}}, \p_{\{1,2,3\}}$ are not stable points by Theorem \ref{thm:plus_one}.

%%%%%%%%%%%%%%%%%%%%%%%%%%%%%%%%%%%%%

\begin{thm} \label{thm:plus_one}
Suppose $(A,\rv)$ is non-degenerate, given $N \times N$ matrix $A$ and $\rv \in \R^N$. 
Assume that $S \subset [N]$ is non-empty and $\p_S \ge \0$ is an asymptotically stable point for Eq \eqref{eq:LV_ODE}. 
If $S \subset T \subseteq [N]$ so that $|T| = |S| + 1$, then $\p_T \ge \0$ is not a stable fixed point for Eq \eqref{eq:LV_ODE}. 
\end{thm}

\begin{proof}
%By Theorem \ref{thm:FT_stable_crit} it can be 
Assume that $S=[N-1]$, without loss of generality. 
It only needs to be shown that $\p$ is not a stable fixed point. 
Since $\p_S$ is stable it is a saturated fixed point so 
\begin{equation}
\begin{aligned}
\e_N^* \p A/{A^S} &= \e_N^* A(\p - \p_S), &\text{ by Proposition \ref{prop:comp_form}}, \\
&< 0, &\text{ by Theorem \ref{thm:FT_stable_crit}}. 
\end{aligned} 
\end{equation}

%%%%%%%%%%%%%%%%%%%%%%%%%%%%%%%%%%%%%

Suppose that $\p$ was also stable. 
Then $\diag(\p)A$ and $\diag(\p^S)A^S$ have positive determinate.
Now by the Schur determinate identity and since $A/{A^S} \in \R$, we have that \\
\begin{equation}
\begin{aligned}
0 &< \det(\diag(\p)A)/\det(\diag(\p^S){A^S}) \\
&= \frac{\det(\diag(\p)) \det(A)}{\det(\diag(\p^S)) \det(A^S)} \\
%= \frac{\det(\diag(\p))}{\det(\diag(\p^S))} \frac{\det(A)}{\det(A^S)} 
&= \e_N^* \p \det(A/ A^S), \text{ Schur determinate identity}, \\
&= \e_N^* \p A/ A^S. 
\end{aligned}
\end{equation}
But then $\e_N^* \p A/{A^S} > 0$, a contradiction. 
Thus $\p$ is not stable. 
\end{proof}

Example \ref{xpl:counter} shows that Theorem \ref{thm:plus_one} doesn't generalize to $|T \setminus S| > 1$. 
Theorem \ref{thm:gen_maxmin} illustrates that the result may be extended in other ways. 

\begin{figure}\centering
\begin{subfigure}{.42\textwidth}\centering
\includegraphics[width=1\linewidth]{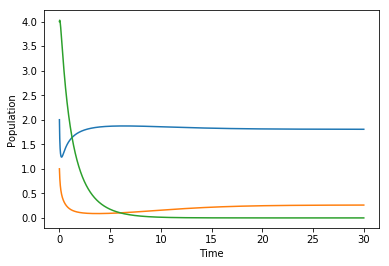}
\end{subfigure}%
\caption{
LV$(A_2,\rv_2)$ has a stable point at $(1.805,0.263,0)$, in Ex. \ref{xpl:counter}. 
} 
\label{fig:contrast_second}
\end{figure}

\begin{xmpl} \label{xpl:counter}

\begin{figure}\centering
\begin{subfigure}{.42\textwidth}\centering
\includegraphics[width=1\linewidth]{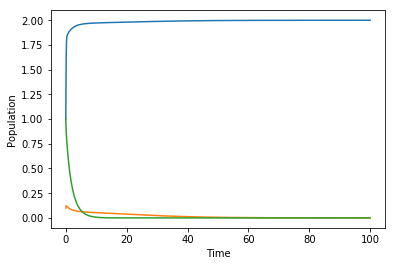}
\end{subfigure}%
\begin{subfigure}{.42\textwidth}\centering
\includegraphics[width=1\linewidth]{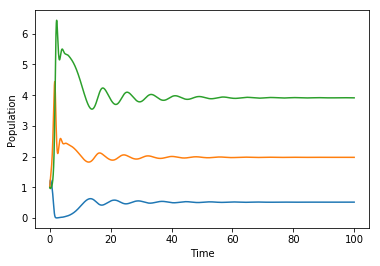}
\end{subfigure}%
\caption{
LV$(A_1,\rv_1)$ has stable points at $(2,0,0)$ and $(0.520833...,1.979166...,3.9166...)$, in Ex. \ref{xpl:counter}. 
} 
\label{fig:counter}
\end{figure}

Consider the two similar Lotka--Volterra systems with $3$ species, LV$(A_1,\rv_1)$ and LV$(A_2,\rv_2)$ where 
\begin{align*}
A_1 =
\begin{pmatrix}
8	& 	4 	& 	1  \\
4.04	& 	1 	& 	1  \\
-1	& 	-1 	& 	0  \\
\end{pmatrix}
\text{ and }
\rv_1 = 
\begin{pmatrix}
16 \\ 
8 \\ 
-2.5 \\ 
\end{pmatrix}
\end{align*}
and
\begin{align*}
A_2 = 
\begin{pmatrix}
  8 & 4.04 & 1 \\
  4.04 & 3 & 1 \\
  -1 & -1 & -0.001 \\
\end{pmatrix},
\rv_2 = 
\begin{pmatrix}
15.5 \\
8.08 \\
-2.5 \\
\end{pmatrix}.
\end{align*}

This is an example, by (Hutson and Vickers, 1983) \cite{hutson1983criterion}, where $S \subset T$ and $|T \setminus S| = 2$, and in fact both $\p_S$ and $\p_T$ are stable. 
The system LV$(A_1,\rv_1)$ has two stable points at $(2,0,0)$ and $(0.520833...,1.979166...,3.9166...)$. 
So clearly Theorem \ref{thm:plus_one} doesn't generalize to $|T \setminus S| > 1$.

There are exactly $2$ stable fixed points for LV$(A_1,\rv_1)$ and the sets of surviving species at these stable points are $\{1\}$ and $\{1,2,3\}$. 
For system LV$(A_2,\rv_2)$, the sets of surviving species at this stable point included $\{1,2\}$ and so by Theorem \ref{thm:plus_one} neither $\{1\}$ nor $\{1,2,3\}$ can be sets of surviving species. 
There is only $1$ stable fixed point for LV$(A_2,\rv_2)$ and the sets of surviving species at this stable point is $\{1,2\}$. 
Although similar, the change from system LV$(A_1,\rv_1)$ to LV$(A_2,\rv_2)$ was enough to alternate the stability each of these fixed points associated with these sets. 
The plots in Figs. \ref{fig:counter} and \ref{fig:contrast_second} and have similar initial conditions, for the purpose of contrasting. 
\end{xmpl}

%%%

The following example applies Theorem \ref{thm:FT_stable_crit} and Theorem \ref{thm:plus_one}. 

\begin{xmpl} \label{xpl:inspr}
Consider the non-degenerate pair $(A, \rv)$ 
\begin{align*}
A =
\begin{pmatrix}
2	& 	1 	& 	1  \\
4	& 	3 	& 	4  \\
2	& 	4 	& 	3  \\
\end{pmatrix}
\text{ and }
\rv = 
\begin{pmatrix}
2 \\ 
5 \\ 
4 \\ 
\end{pmatrix}.
\end{align*}

This example is interesting because the points $\p_{\{1\}},\p_{\{1,2\}},\p_{\{1,2,3\}} \in \R_{>0}^N$ with sets $\{1\}\subset \{1,2\} \subset \{1,2,3\}$ with cardinality one apart. 
The point $\p_{\{1,2\}} = (1/2,1,0)^T$ and so $\diag(\p^{\{1,2\}}) A^{\{1,2\}}= 
\begin{pmatrix} 
1	& 	1/2 	\\
4	& 	3 	\\
\end{pmatrix}$ has all of its eigenvalues with positive real part, and $A(\p - \p_{\{1,2\}})=(0,0,-1)^T \le \0$. 
Thus $\p_{\{1,2\}}$ is asymptotically stable. 

The point $\p = \p_{\{1,2,3\}}$ is a saturated fixed point since $A(\p - \p_{\{1,2,3\}})= \0$, 
and the point $\p_{\{1\}}=(1,0,0)^T$ has that $-\diag(\p^{\{1\}}) A^{\{1\}}=-2$ is stable, 
but $\p$ and $\p_{\{1\}}$ are not stable by Theorem \ref{thm:plus_one}. 
\end{xmpl}

%%%%%%%%%%%%%%%%%%%%%%%%%%%%%%%%%%%%%

%%%%%%%%%%%%%%%%%%%%%%%%%%%%%%%%%%%%%
%%%%%%%%%%%%%%%%%%%%%%%%%%%%%%%%%%%%%
%\pagebreak

%%%%%%%%%%%%%%%%%%%%%%%%%%%%%%%%%%%%%
%%%%%%%%%%%%%%%%%%%%%%%%%%%%%%%%%%%%%

%Subsection at end or appendix 

Consider the following assumption, given the $N \times N$ matrix $A$, the vector $\rv \in \R^N$, and $S \subset T \subseteq [N]$: 
\bi
\item \textbf{Assumption A}$_{S,T}$:  Suppose $(A, \rv)$ is non-degenerate and that $S, T$ are non-empty so that 
\begin{align*}
\l(A^{S,T \setminus S} / A^S \r) \l(P^{T \setminus S}\r)^* \p_T \notin \R_{<0}^{|T \setminus S|}. 
\end{align*}
\ei
Theorem \ref{thm:gen_maxmin} illustrates how \textbf{A}$_{S,T}$ might be used to discover new properties of stability of Lotka--Volterra.

%%%% 

\begin{thm} \label{thm:gen_maxmin}
Suppose \textbf{A}$_{S,T}$ of $(A, \rv)$. 
If $S \subset T \subseteq [N]$ and $\p_S$ is a stable point of LV$(A, \rv)$ then $\p_T$ is not a stable point of LV$(A, \rv)$. 
In fact, $\p^T$ is not a stable point of LV$(A^T, \rv^T)$. 
\end{thm}

\begin{proof}
Suppose \textbf{A}$_{S,T}$ of $(A, \rv)$ with $S \subset T \subseteq [N]$ so that $\p_S > \0$ is a stable point of LV$(A, \rv)$. 
By Theorem \ref{thm:FT_stable_crit} it is enough to show that $\p^T > \0$ is not a stable point of LV$(A^T, \rv^T)$. 

Suppose $\p^T$ is a stable point of LV$(A^T, \rv^T)$. 
Then $\l(P^{T \setminus S}\r)^* A_T (\p_T - \p_S) \le \0$ and since $(A, \rv)$ is non-degenerate, $\l(P^{T \setminus S}\r)^* A (\p_T - \p_S) < \0$. 
Thus by Proposition \ref{prop:comp_form} 
\begin{align*}
\l(A^{S,T \setminus S} / A^S \r) \l(P^{T \setminus S}\r)^* \p_T 
&= \l(P^{T \setminus S}\r)^* A (\p_T - \p_S) < \0,
\end{align*}
contradicting \textbf{A}$_{S,T}$. 
Thus $\p^T$ is not a stable point of LV$(A^T, \rv^T)$. 
\end{proof}

%%%%

In Example \ref{ex:rand_A} the set of all stable points are given. 
The sets of surviving species at the stable points do not contain any other in this example, which is common. 
The values of the community matrix are chosen to be somewhat evenly distributed. 
Theorem \ref{thm:gen_maxmin} is used to show that there $\p$ is not stable.

\begin{figure}\centering
\begin{subfigure}{.33\textwidth}\centering
\includegraphics[width=1\linewidth]{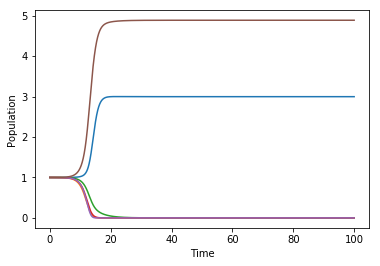}
\end{subfigure}%
\begin{subfigure}{.33\textwidth}\centering
\includegraphics[width=1\linewidth]{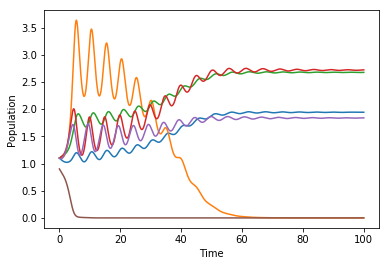}
\end{subfigure}%
\begin{subfigure}{.33\textwidth}\centering
\includegraphics[width=1\linewidth]{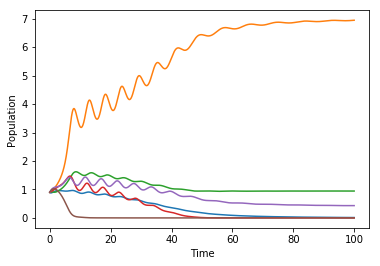}
\end{subfigure}%
\caption{
The Lotka--Volterra system given by $A$ and $\p$ in Ex. \ref{ex:rand_A}. 
} 
\label{fig:example_criterion}
\end{figure}

\begin{xmpl} \label{ex:rand_A}
Suppose there are $6$ species in the Lotka--Volterra system with matrix $A$ and vector $\rv$ so that  
\begin{align*}
A = 
\begin{pmatrix}
  0.99 & -0.04 & -0.846 & 0.069 & -0.677 & -0.893\\
  0.518 & 0.025 & 0.716 & -0.51 & -0.316 & 0.288\\
  -0.544 & -0.014 & 0.822 & -0.375 & 0.39 & 0.561\\
  -0.144 & 0.17 & 0.472 & 0.45 & -0.925 & 0.484\\
  0.808 & 0.279 & 0.268 & -0.164 & 0.255 & 0.863\\
  -0.353 & 0.26 & 0.593 & 0.937 & 0.308 & 0.721\\
\end{pmatrix},
\p = 
\begin{pmatrix}
1 \\
1 \\
1 \\
1 \\
1 \\
1 \\
\end{pmatrix}.
\end{align*}

\noindent

The three plots in Fig. \ref{fig:example_criterion} show the evolution of the dynamical system starting near the fixed point $\p$. 
The solutions converge to a stable fixed point after some species has died off. 
From left to right the surviving species are $\{1, 6\}$, $\{2, 3, 5\},$ and $\{1, 3, 4, 5\}$, none of which contain another. 
The periodic waves indicate that the $6$ dimensional population vectors are rotating, essentially spiraling towards the stable fixed point. 

Now, Theorem \ref{thm:gen_maxmin} can be used to show that $\p$ is unstable. 
Consider the Schur complement 
\begin{align*}
A^{S,S^c}/A^S = 
\begin{pmatrix}
  0.632 & -1.014 & -2.54 \\
  0.015 & -0.412 & 1.317 \\
  -1.085 & -1.345 & -1.153 \\
\end{pmatrix}
\end{align*}

and so 
\begin{align*}
\l(A^{S,S^c} / A^S \r) \l(P^{S^c}\r)^* \p = 
\begin{pmatrix}
-2.922 \\
0.92 \\
-3.583 \\
\end{pmatrix} 
\notin \R_{<0}^{|S^c|}. 
\end{align*}
So by Theorem \ref{thm:gen_maxmin} $\p$ is not stable, because $\p_S$ is stable. 
\end{xmpl}

%%%%

The following results illustrate another case where \textbf{A}$_{S,T}$ applies. 

%%%%

\begin{lem} \label{lem:symequiv_posets}
Assume $(A, \rv)$ is non-degenerate and that $A$ is symmetric. 
Suppose that $S \subseteq [N]$ and $\p_S \ge 0$. 
Then $A^S > 0$ if and only if $-\diag(\p^S) A^S$ is stable. 
\end{lem}

The proof is given immediate because stable symmetric matrices are $D$-stable, see Lemma 3.2.1 in \cite{takeuchi1996global}. 

\iffalse
\begin{proof}
Fix $S \subseteq [N]$ to be nonempty, so that $\p_S \ge \0$. 
Suppose that $A^S$ and $\diag(\p^S)$ are symmetric and $\diag(\p^S)$ 
is positive definite. 
It is well known that the $k$th eigenvalue 
$\lam_k(\diag(\p^S) A^S) = \lam_k \l(\sqrt{\diag(\p^S)} A^S \sqrt{\diag(\p^S)}\r)$. 
It is a simple application of the Courant--Fischer min-max principle and adjoint properties of an inner product. 
Hence $-\diag(\p^S) A^S$ is stable if and only if $-\sqrt{\diag(\p^S)} A^S \sqrt{\diag(\p^S)}$ is stable. 
Another well known fact is that if $A^S$ and $\diag(\p^S)$ are symmetric and $\diag(\p^S)$ is positive definite then 
$\sqrt{\diag(\p^S)} A^S \sqrt{\diag(\p^S)} > 0$ if and only if $A^S > 0$. 
Thus $A^S > 0$ if and only if $-\diag(\p^S) A^S$ is stable. 
\end{proof}
\fi

%%%%

\begin{lem} \label{lem:sym_inviolable}
Suppose $(A, \rv)$ is non-degenerate for symmetric $A>0$, then for any $S \subset T \subseteq [N]$, $(A, \rv)$ satisfies \textbf{A}$_{S,T}$. 
\end{lem}

\begin{proof}
Suppose that there is a non-empty $S \subset T \subset [N]$ such that 
$-\diag(\p^S) A^S$ is stable and $\p_S \ge \0$. 
By Lemma \ref{lem:symequiv_posets} $A^S, A^T > 0$. 

By Haynsworth's inertia additivity formula \ref{prop:inertia}, 
\begin{align}
\In\l(A^{S,T \setminus S} / A^S \r) = (|S^c|, 0, 0), 
\end{align}
thus the matrix is positive definite. 
This means that 
\begin{align}
\diag \l(\l(P^{T \setminus S}\r)^* \p_T \r) \l(A^{S,T \setminus S} / A^S \r) 
\diag \l(\l(P^{T \setminus S}\r)^* \p_T \r)
\end{align}
is positive definite as well. 
By definition of positive definite 
\begin{align}
\1^* \diag \l(\l(P^{T \setminus S}\r)^* \p_T \r) \l(A^{S,T \setminus S} / A^S \r) 
\diag \l(\l(P^{T \setminus S}\r)^* \p_T \r) \1 > 0. 
\end{align}
Thus 
\begin{align}
\diag \l(\l(P^{T \setminus S}\r)^* \p_T \r) \l(A^{S,T \setminus S} / A^S \r) 
\diag \l(\l(P^{T \setminus S}\r)^* \p_T \r) \1 < \0,
\end{align}
is false. 
Hence the assumption \textbf{A}$_{S,T}$ is satisfied, else $\p_T$ is not in the positive orthant. 
\end{proof}

%%%%

%%%%

\begin{cor} \label{cor:FT_sym_crit}
Suppose $(A, \rv)$ is non-degenerate. 
If $S \subset T \subseteq [N]$, $A_T$ is symmetric, and $\p_S$ is a stable point of LV$(A, \rv)$, then $\p_T$ is not a stable point of LV$(A, \rv)$. 
\end{cor}
\begin{proof}
Immediate from Lemmas \ref{lem:symequiv_posets} and \ref{lem:sym_inviolable}, and Theorem \ref{thm:gen_maxmin}. 
\end{proof}

% This result was already shown in section \ref{ch:ring}, Theorem \ref{thm:stable_sub}. 
Corollary \ref{cor:FT_sym_crit} was completely generalized in work regarding P-matrices in \cite{takeuchi1996global}, which depended on Murty's Theorem \ref{thm:Murty_P}. 
However, Theorem \ref{thm:gen_maxmin} does not depend on Murty's Theorem.

\iffalse %TURN ON
\begin{thm} \label{thm:comb_sym}
Suppose that $A$ is irreducible and has only cycles of length $\ge 2$. 
Then system LV$(A,\rv)$ and every {\em reduced system} have a nonnegative and globally stable equilibrium point for each $\rv \in \R^N$ iff $-A$ is a $P$-matrix. 
\end{thm}
\fi

%%%%%%%%%%%%%%%%%%%%%%%%%%%%%%%%%%%%%
%%%%%%%%%%%%%%%%%%%%%%%%%%%%%%%%%%%%%

\section{Conclusion} \label{sec:conc}

Relationships within the community matrix were shown to predict stability for the generalized Lotka--Volterra system. 
A strong relationship exists between the transversal eigenvalues with respect to a subset of the species in a system and the corresponding Schur complement of the Jacobian at the fixed point with the submatrix determined by the same subset of species. 
This relationship was exploited throughout giving rise to new methods to check the stability of fixed points.

%%%%%%%%%%%%%%%%%%%%%%%%%%%%%%%%%%%%%
%%%%%%%%%%%%%%%%%%%%%%%%%%%%%%%%%%%%%
\pagebreak

\appendix

\section{Appendix} \label{apx:LV}

%\section{Prerequisites of Lotka--Volterra} \label{apx:LV}
%\section*{} \label{apx:LV}

%%%%%%%%%%%%%%%%%%%%%%%%%%%%%%%%%%%%%
%%%%%%%%%%%%%%%%%%%%%%%%%%%%%%%%%%%%%

A matrix $P$ is a \textit{projection} matrix if $P^2 = P$. 
Furthermore, $P$ is an \textit{orthogonal projection} if $P^2 = P = P^*$. 
A matrix $M$ is said to be a $P$-matrix if it is a complex valued square matrix with every principal minor greater than 0. 
Finally, $M$ is said to be \textit{stable} if the real part of all of its eigenvalues are negative. 

Given $N$ species, a solution to \eqref{eq:LV_ODE}, $\x(t)$, will remain in $\R_{>0}^N$. 
Yet it is possible for $x_i(t)$ to get arbitrarily close to zero for $t \in \R_{\ge 0}$, and move away again. 
To that end consider the following definitions. 

\begin{define} \label{def:permanent}
The system given by \eqref{eq:LV_ODE} is {\em permanent} if there exists an $\epsilon >0$ so that for any $\x(0) > \0$, $\liminf_{t \to \infty} \x(t) > \epsilon \1$ and $\limsup_{t \to \infty} \x(t) > \frac{1}{\epsilon} \1$. 
\end{define}

\begin{define} \label{def:persistent}
The system given by \eqref{eq:LV_ODE} is {\em persistent} if there exists an $\epsilon >0$ so that for any $\x(0) > \0$ and $i \in [N]$, $\limsup_{t \to \infty} x_i(t) > \epsilon$. 
\end{define}

A persistent system means that the population is compactly bounded above zero, but this bound may depend on the initial condition. 
A permanent system means that each of its populations are bounded away from zero. 
Clearly, permanent implies persistent.

\begin{define} \label{def:comp_set}
Suppose that $A$ is a real valued $N \times N$ matrix, and $\rv \in \R^N$. 
A fixed point $\x$ is a {\em saturated} fixed point for the Lotka--Volterra system LV$(A, \rv)$ if 
$\rv - A \x \le \0$.
\end{define}

\begin{define} \label{def:pos_def}
An $N \times N$ real valued matrix $A$ (many authors require $A$ to be symmetric) is positive definite if 
\begin{align*}
\ip{A\x}{\x} > 0, \forall \x \in \R^N,
\end{align*}
and positive semi-definite if 
\begin{align*}
\ip{A\x}{\x} \ge 0, \forall \x \in \R^N.
\end{align*}
Denote a positive definite matrix by $A > 0$ and a positive semi-definite matrix by $A \ge 0$. 
For $N \times N$ matrices $A$ and $B$, $A > B$ is written to mean that $A - B > 0$ (respectively for positive semi-definite).
\end{define}

\begin{define} \label{def:SchurCompl}
Suppose the matrices $A,B,C,D$ have dimensions $n \times n$, $n \times m$, $m \times n$, $m \times m$ respectively. 
Let $M$ be the $n+m$ by $n+m$ block matrix 
\begin{align*}
M = 
\begin{pmatrix}
A & B \\
C & D \\
\end{pmatrix}.
\end{align*} 
Define the {\em Schur complement} to be the block $D$ of the matrix $M$ is the $m \times m$ matrix denoted 
$M/D = A - BD^{\dagger}C$. 
Likewise the Schur complement of the block $A$ of the matrix $M$ is the $n \times n$ matrix denoted 
$M/A = D - CA^{\dagger}B$. 
\end{define}

%%%%%%%%%%%%%%%%%%%%%%%%%%%%%%%%%%%%%
%%%%%%%%%%%%%%%%%%%%%%%%%%%%%%%%%%%%%

\begin{define} \label{def:LCP}
The linear complementarity problem (sometimes denoted as LCP$(A, \rv)$) is the system of inequalities 
\begin{align} \label{eq:LCP}
A \x - \rv &\ge \0 \nonumber \\
\x &\ge \0 \nonumber \\
\ip{A \x - \rv}{\x} &= 0,
\end{align}
given $A$ is an $N \times N$ matrix and $\rv \in \R^N$. 
\end{define}
LCPs are often used to model contact forces between rigid bodies, and are used in many applied industrial problems. 
There need not always be a solution $\x$ to LCP$(A, \rv)$. 

\begin{prop} \label{prop:LV_LCP}
Solutions to LCP$(A,\rv)$ are equivalent to saturated equilibrium points for the LV$(A,\rv)$, \cite{takeuchi1996global}. 
\end{prop}

The proof of Theorem \ref{prop:LV_LCP} is in \cite{takeuchi1996global} and is immediate from Theorem \ref{thm:stable_imply_saturated}. 

\begin{thm} \label{thm:stable_imply_saturated}
Consider the system
$\dot x_i = x_i f_i(x_1,...,x_n), \forall i \in [n]$
where $f_i$'s are continuous. 
If $\p \ge \0$ is stable, then $f_i(p_1,...,p_n) \le 0, \forall i \in [n]$.
\end{thm}

Murty's Theorem \ref{thm:Murty_P} \cite{murty1988linear} proven by Murty in 1972 proves to be a powerful tool for finding stable points for the Lotka--Volterra. 

\begin{thm} \label{thm:Murty_P}
The LCP$(A,\rv)$ has a unique solution for each $\rv \in \R^N$ iff $A$ is a $P$-matrix. 
\end{thm}

\pagebreak

%%%%%%%%%%%%%%%%%%%%%%%%%%%%%%%%%%%%%
%%%%%%%%%%%%%%%%%%%%%%%%%%%%%%%%%%%%%

%\section{Exploiting Jacobian of Lotka--Volterra} \label{ch:ap_2}

\section{Appendix} \label{apx:props}

Some necessary properties are covered in this appendix. 

\begin{prop}
Given that $(A,\rv)$ is non-degenerate, the following are true:

\begin{enumerate}
	\item \label{1} $\p_S = P^S \p^S$
	\item \label{2} $(P^T)^* A P_T = A^T (P^T)^*$
	\item \label{3} $(P^T)^* A(\p - \p_T) = \0$
	\item \label{4} $(P^S)^* A_{S \cup \{ i \} } \e_i = (P^S)^* A \e_i $
	\item \label{5} $A_S^{\dagger} = P^S((P^S)^* A P^S)^{-1} (P^S)^*$ 
	\item \label{6} $(A_S^{\dagger})_{S, T} = \begin{pmatrix}
(A^S)^{-1} 		& 0_{|S| \times |T|} \\
0_{|T| \times |S|}		& 0_{|T| \times |T|} \\
\end{pmatrix}$
	\item \label{7} $A^{T} (P^T)^* \p_S = A^{T} (P^T)^* A_S^{\dagger} P^T A^{T} \p^T$, for $S \subset T$
	\item \label{8} $A^{S,T} - A^{S,T} (A_S^{\dagger})_{S, T} A^{S,T} = 
\begin{pmatrix}
0_{|S| \times |S|} 	& 0_{|S| \times |T|} \\
0_{|T| \times |S|} 	&  A^{S,T} / A^{S} \\
\end{pmatrix}$, for disjoint $S$ and $T$. 
\end{enumerate}
\end{prop}

\begin{proof}
The first six are trivial and left to the reader. 
For item \ref{7} 
\begin{align*} 
A^{T} (P^T)^* \p_S &= A^{T} (P^T)^* A_S^{\dagger} \rv_S \\
&= A^{T} (P^T)^* A_S^{\dagger} \rv_T \\
&= A^{T} (P^T)^* A_S^{\dagger} P^T A^{T} \p^T .
\end{align*}
For item \ref{8} 
\begin{align*} 
A^{S,T} - &A^{S,T} (A_S^{\dagger})_{S, T} A^{S,T} 
= 
A^{S,T}
-
A^{S,T}
\begin{pmatrix}
(A^S)^{-1} 		& 0_{|S| \times |T|} \\
0_{|T| \times |S|}		& 0_{|T| \times |T|} \\
\end{pmatrix} 
A^{S,T}, \mbox{ by item } \ref{6}, \\
&= 
\begin{pmatrix}
(P^S)^* A P^S 	&	(P^S)^* A P^T \\
(P^T)^* A P^S	& 	(P^T)^*A P^T \\
\end{pmatrix} 
-
\begin{pmatrix}
(P^S)^* A P^S 		&	(P^S)^* A P^T \\
(P^T)^* A P^S	& 	(P^T)^* A P^S (A^S)^{-1} (P^S)^* A P^T  \\
\end{pmatrix}, \mbox{ product, } \\
& = 
\begin{pmatrix}
0_{|S| \times |S|} 	& 0_{|S| \times |T|} \\
0_{|T| \times |S|} 	& A^{S,T} / A^{S} \\
\end{pmatrix}, \mbox{ definition}.
\end{align*}
\end{proof}

The following two propositions are exploited in section $5$. 
As can be seen, the first of these propositions relates to the property of saturated fixed points, 
and the second is a property of the Schur complement of two Jacobians. 

%\pagebreak

\begin{prop} \label{prop:comp_form}
Suppose $(A,\rv)$ is non-degenerate. 
If $S \subset T \subseteq [N]$ is non-empty then 
\begin{align*}
(P^{T \setminus S})^* A (\p_T - \p_S) = (A^{S, T \setminus S} / A^S) (P^{T \setminus S})^* \p_T. 
\end{align*}
\end{prop}

\begin{proof}
Using the above properties gives us

\begin{align*}
(P^{T \setminus S})^* A (\p_T - \p_S) 
&= (P^{T \setminus S})^* P_T A P_T (\p_T - \p_S), \mbox{ by item } \ref{1}, \\
&= (P^{T \setminus S})^* P^T A^T (P^T)^* (\p_T - \p_S), \mbox{ by item } \ref{2}, \\
&= (P^{T \setminus S})^* P^T (A^T - A^T (P^T)^* A_S^{\dagger} P^T A^T) \p^T, \mbox{ by item } \ref{7}, \\
&= (P^{T \setminus S})^* P^T (A^T - A^T (P^T)^* A_S^{\dagger} P^T A^T) (P^T)^* \p_T \\
&= (P^{T \setminus S})^* (A_T - A_T A_S^{\dagger} A_T)  \p_T \\
&= (P^{T \setminus S})^* P_T (A - A P_T A_S^{\dagger} P_T A) P_T  \p_T \\
&= (P^{T \setminus S})^* P^{S,T \setminus S}(A^{S,T \setminus S} - A^{S,T \setminus S} (A_S^{\dagger})_{S, T \setminus S} A^{S,T \setminus S}) (P^{S,T \setminus S})^* \p_T \\
&= (P^{T \setminus S})^* P^{S,T \setminus S} 
\begin{pmatrix}
0_{|S| \times |S|} 			& 0_{|S| \times |T \setminus S|} \\
0_{|T \setminus S| \times |S|} 	&  A^{S,T} / A^{S} \\
\end{pmatrix}
(P^{S,T \setminus S})^* \p_T, \mbox{ by item } \ref{8}, \\
&= 
\begin{pmatrix}
0_{|T \setminus S| \times |S|} 	&  A^{S,T \setminus S} / A^{S} \\
\end{pmatrix}
\begin{pmatrix}
(P^S)^* \p_T \\
(P^{T \setminus S})^* \p_T \\
\end{pmatrix} \\
&= (A^{S, T \setminus S} / A^S) (P^{T \setminus S})^* \p_T. 
\end{align*}
\end{proof}

\begin{prop} \label{prop:stable_form}
Suppose $(A,\rv)$ is non-degenerate. 
If $S \subset T \subseteq [N]$ is non-empty then 
\begin{align*}
(\diag(\p_T) A_T)^{S,T \setminus S} / (\diag(\p_T) A_T)^S = \diag((P^{T \setminus S})^* \p_T) (A^{S, T \setminus S} / A^S). 
\end{align*}
\end{prop}

\begin{proof}

First notice that 
\begin{align*}
(\diag(\p_T) A_T)^{S,T \setminus S} 
&= (P^{S,T \setminus S})^* \diag(\p_T) A_T P^{S,T \setminus S} \\
&= (P^{S,T \setminus S})^* \diag(\p_T) P^{S,T \setminus S} (P^{S,T \setminus S})^* A_T P^{S,T \setminus S} \\
&= (\diag(\p_T))^{S,T \setminus S} A^{S,T \setminus S} \\
&= 
\begin{pmatrix}
\diag((P^S)^* \p_T) 				& 0 \\
0 							& \diag((P^{T \setminus S})^* \p_T) \\
\end{pmatrix}
\begin{pmatrix}
A^S 		&	(P^S)^* A P^{T \setminus S} \\
(P^{T \setminus S})^* A P^S	& 	A^{T \setminus S} \\
\end{pmatrix} \\
&= 
\begin{pmatrix}
\diag((P^S)^* \p_T) A^S 						&	\diag((P^S)^* \p_T) (P^S)^* A P^{T \setminus S} \\
\diag((P^{T \setminus S})^* \p_T) (P^{T \setminus S})^* A P^S	& 	\diag((P^{T \setminus S})^* \p_T) A^{T \setminus S} \\
\end{pmatrix}.
\end{align*}
Thus 
\begin{align*}
(\diag(\p_T) A_T)^{S,T \setminus S} &/ (\diag(\p_T) A_T)^S = \\
&= \diag((P^{T \setminus S})^* \p_T) A^{T \setminus S} - ...\\
- \diag((P^{T \setminus S})^* \p_T) (P^{T \setminus S})^* &A P^S (\diag((P^S)^* \p_T) A^S )^{-1} \diag((P^S)^* \p_T) (P^S)^* A P^{T \setminus S} \\
&= \diag((P^{T \setminus S})^* \p_T) A^{T \setminus S} - \diag((P^{T \setminus S})^* \p_T) (P^{T \setminus S})^* A P^S (A^S )^{-1} (P^S)^* A P^{T \setminus S} \\
&= \diag((P^{T \setminus S})^* \p_T) (A^{S, T \setminus S} / A^S).
\end{align*}
\end{proof}

%\pagebreak

Due to the use of principal submatrices, the following results are relevant. 

\begin{prop} \label{prop:submat}
Suppose $S \subset T \subseteq [N]$ then $A^S$ is a principal submatrix of $A^T$. 
\end{prop}

\begin{proof}
Notice that $P^S P_T = P^S$, so $A^S = P^S A (P^S)^* = P^S P_T A P_T^* (P^S)^* = P^S (P^T)^* A^T P^T (P^S)^*$.
Thus $A^S$ is a $|S| \times |S|$ principal submatrix of $A^T$.
\end{proof}

\begin{prop} \label{prop:incl}
Suppose $S \subset T\subseteq [N]$ and $A^{T}$ is a $P$-matrix, then $A^{S}$ is a $P$-matrix (respectively symmetric positive definite).
\end{prop}

\begin{proof}
By Proposition \ref{prop:submat} and the definition of $P$-matrix (respectively symmetric positive definite).
\end{proof}

The following is the Schur determinant identity \cite{crabtree1969identity}. 
\begin{prop} \label{prop:Schur_det_id}
Consider the block matrix $M$ and assume $A$ is nonsingular,
\begin{align*}
M = 
\begin{pmatrix}
A & B \\
C & D \\
\end{pmatrix}.
\end{align*}
Then $\det(M) = \det(A) \det(D - C A^{-1}B)$.
\end{prop}

Define the inertia of a matrix $In(M) := (\pi(M),\nu(M),\delta(M))$ where $\pi,\nu,\delta$ are the number of positive, negative and zero eigenvalues of the input matrix respectively. 
The following is known as the Haynsworth's inertia additivity formula \cite{puntanen2005historical},
\begin{prop} \label{prop:inertia}
Let 
\begin{align*}
M = 
\begin{pmatrix}
A & B \\
B^* & D \\
\end{pmatrix}
\end{align*} 
 be a self-adjoint matrix. 
Then the inertia of $M$ is 
\begin{align*}
\In(M) = \In(A) + \In(M/A) = \In(A) + \In(D - B^* A^{-1} B). 
\end{align*}
\end{prop}

\begin{prop} \label{prop:comp_nondeg}
Suppose $A$ is non-degenerate and $S \subset [N]$, then $\l(A_{ii} - \e_i^* A A_S^{\dagger}  A \e_i \r) \neq 0$, for each $i \in S^c $. 
\end{prop}

\begin{proof}
Let $i \in S^c $. 
Then 
\begin{align*}
A^{S,i} &:= 
\begin{pmatrix}
A^S 			&	(P^S)^* A_S \e_i \\
\e_i^* A_S P^S 	& 	A_{ii} \\
\end{pmatrix}.
\end{align*}
So by Proposition \ref{prop:Schur_det_id} 
\begin{align*}
\det(A^{S \cup \{i\}}) = \det(A^{S,i}) = \det(A^S) \l(A_{ii} - \e_i^* A A_S^{\dagger}  A \e_i \r).
\end{align*}
By definition both $\det(A^{S \cup \{i\}})$ and $\det(A^S)$ are nonzero, hence the conclusion. 
\end{proof}

%%%%%%%%%%%%%%%%%%%%%%%%%%%%%%%%%%%%%
%%%%%%%%%%%%%%%%%%%%%%%%%%%%%%%%%%%%%

\pagebreak

\bibliography{thesis_bib}
\bibliographystyle{amsalpha}

\end{document}